\documentclass[12pt]{article}
\usepackage{amssymb,amsfonts}
\usepackage{color}
\usepackage{amsfonts}
\usepackage{latexsym}
\usepackage[dvips]{graphics}
\usepackage{epsfig}
\usepackage{mathdots}
\usepackage{bm}
\usepackage{amsmath}
\usepackage{mathrsfs}
\usepackage{lineno}
\newtheorem{lem}{Lemma}
\newtheorem{theo}{Theorem}
\newtheorem{pro}{Proposition}
\newtheorem{cor}{Corollary}
\newtheorem{cla}{Claim}

\newcommand{\proof}{{\noindent {\em Proof}.\quad}\setcounter{countclaim}{0}
\setcounter{countcase}{0}}
\newcommand{\proofend}{{\hfill$\Box$}}

\newcounter{countfig}

\newcounter{countclaim}

\newcounter{countcase}

\newcommand{\beeq}{\begin{equation}}
\newcommand{\eneq}{\end{equation}}

\newcommand{\beeqn}{\begin{eqnarray*}}
\newcommand{\eneqn}{\end{eqnarray*}}

\newcommand {\red} {\textcolor{red}}
\newcommand {\green} {\textcolor{green}}
\newcommand {\blue} {\textcolor{blue}}

\def \iff {if and only if }

\def \usecolour
{\newcommand {\rered}{\red}
\newcommand {\reblue} {\blue}
\newcommand {\regreen} {\green}
}

\usecolour   
\newcommand {\relabel}[1] {\label{#1} \red{[*: #1]}}

\def\relabel {\label}   

\def \M {\mathscr{M}}

\begin{document}

\newcommand{\resection}[1]
{\section{#1}\setcounter{equation}{0}}

\renewcommand{\theequation}{\thesection.\arabic{equation}}

\baselineskip 0.6 cm

\title {\bf New upper bounds for the crossing numbers of crossing-critical graphs\thanks {The work was supported by the National
Natural Science Foundation of China (No. 11301169),
Hunan Provincial Natural Science Foundation of China (No.
2017JJ2055), and the Scientific Research Fund of Hunan Provincial Education Department (No. 19B116).}}

\author
{Zongpeng Ding\thanks{Email: dzppxl@163.com.},
 Zhangdong Ouyang\thanks{Corresponding author. Email: oymath@163.com.}\\
\small School of mathematics and Computational Sciences\\
\small Hunan First Normal University,
Changsha 410205, P.R China\\
\\
Yuanqiu Huang\thanks{Email: hyqq@hunnu.edu.cn.}\\
\small Department of Mathematics\\
\small Hunan Normal University, Changsha 410081, P.R China\\
\\
Fengming Dong\thanks{Email: fengming.dong@nie.edu.sg (expired on 24/03/2027)
and donggraph@163.com.}\\
\small National Institute of Education,  Nanyang Technological
University,  Singapore
}

\date{}

\maketitle

\begin{abstract}
A graph $G$ is {$k$-crossing-critical} if
$cr(G)\ge k$, but
$cr(G\setminus e)<k$ for each edge $e\in E(G)$,
where $cr(G)$ is the crossing number of $G$.
It is known that for any $k$-crossing-critical graph $G$,
$cr(G)\le 2.5k+16$ holds,
and in particular, if $\delta(G)\ge 4$, then
$cr(G)\le 2k+35$ holds,
where $\delta(G)$ is the minimum degree of $G$.
In this paper, we improve these upper bounds to
$2.5k +2.5$ and $2k+8$ respectively.
In particular, for any $k$-crossing-critical graph $G$
with $n$ vertices, if $\delta(G)\ge 5$, then
$cr(G)\le 2k-\sqrt k/2n+35/6$ holds.
\end{abstract}




\resection{Introduction \relabel{sec1}}

All graphs considered here are simple, connected, finite and undirected unless otherwise specified.
For any graph $G$,
let $V(G)$, $E(G)$ and $\delta(G)$ denote its vertex set, edge set and minimum degree.
A {\it drawing} of a graph $G$ is
a mapping $D$ that assigns to each vertex in $V(G)$
a distinct point in the plane,
and to each edge $uv$ in $G$
a continuous arc connecting $D(u)$ and $D(v)$,
not passing through the image of any other vertex.
For any drawing $D$ of $G$,
let $cr(D)$ denote the number of crossings in $D$,
and the {\it crossing number} of $G$,
denoted by $cr(G)$,  is the minimum value of $cr(D)$'s
among all possible drawings $D$ of $G$. For more on crossing numbers of graphs, we refer to \cite{MS} and the references therein.

A graph $G$ is {\it $k$-crossing-critical}
if $cr(G)\ge k$,
but $cr(G\setminus e)<k$ for every edge $e\in E(G)$
(e.g. see \cite{ZD}).
A graph is {\it crossing-critical}
if it is $k$-crossing-critical for some $k$.


Crossing-critical graphs give insight into structural properties of the crossing
number invariant and have thus generated considerable interest.
Let $\M_k$   denote the set of $k$-crossing-critical graphs.
Richter and Thomassen~\cite{RB} showed that
$cr(G)\le 2.5k+16$ holds for each $G\in \M_k$.
Salazar~\cite{GS} improved this result to $cr(G)\le 2k+35$
for the case that $\delta(G)\ge 4$.
Lomel\'i and Salazar~\cite{LS} showed that,
for each integer $k > 0$,
there is an integer $n_k$ such that
for any $G\in \M_k$ with at least $n_k$ vertices of degree greater than two,
$cr(G) \le 2k + 23$ holds.

In this paper, we further improve
Richter and Thomassen's result in \cite{RB} and
Salazar's result in \cite{GS} to
$cr(G)\le 2.5k+2.5$ and $cr(G)\le 2k+8$ respectively.
Furthermore,  we show that,
for any $G\in \M_k$ with $n$ vertices,
if $\delta(G)\ge 5$, then
$cr(G)\le 2k-\sqrt k/2n+35/6$ holds.

\resection{Choose a suitable cycle in a graph \label{sec2}
}

For any positive integers $l$ and $\Delta$,
a cycle $C$ in a graph $G$ is called {\it $(l,\Delta)$-nearly-light}
if the length of $C$ is at most $l$ and
at most one vertex of $C$
has degree greater than $\Delta$
(see \cite{LS}).

For any cycle $C$ in a graph $G$,  define
$$
\mu(C)=\min_{v\in V(C)}
\sum_{u\in V(C)\setminus\{v\}}(d_G(u)-2).
$$
Clearly, if $C$ is an $(l,\Delta)$-nearly-light cycle in $G$,
then $\mu(C)\le (l-1)(\Delta-2)$.
In \cite{LS},
Lomel\'i and Salazar  showed that
if a $k$-crossing-critical graph $G$ contains a cycle $C$
with $\mu(C)\le s$,
then $cr(G)\le 2(k-1)+s/2$, and therefore,
if $G$ contains an $(l,\Delta)$-nearly-light cycle,
then  $cr(G)\le 2(k-1)+(l-1)(\Delta-2)/2$ holds.

The {\it skewness} of a graph $G$, denoted by $sk(G)$,
is the minimum integer $t$ such that
$G\setminus E_0$ is planar for a subset $E_0$ of $E(G)$
with $|E_0|=t$.
In this section, we show that for any graph $G$
with $\delta(G)\ge 3$,
$G$ contains a cycle $C$ with $\mu(C)\le sk(G)+10$.
In Section~\ref{sec4}, we apply this fact to find
an upper bound for $cr(G)$ in terms of $sk(G)$ and $\delta(G)$,
where $G$ is crossing-critical.

The next elementary result will be applied later.

\begin{lem}\relabel{lem1}
Let $d_1,d_2,\cdots,d_5$ be integers with
$3\le d_1\le d_2\le \cdots \le d_5$.
Then, the following statements hold:
\begin{enumerate}
\renewcommand{\theenumi}{\rm (\arabic{enumi})}
\item if $\sum_{i=1}^{3}\frac{1}{d_i}>\frac{1}{2}$, then $\sum_{i=1}^{2}(d_i-2)\le 10$;
\item 
if $\sum_{i=1}^{4}\frac{1}{d_i}>1$, then $\sum_{i=1}^{3}(d_i-2)\le 5$; and
\item 
if $\sum_{i=1}^{5}\frac{1}{d_i}>\frac{3}{2}$, then $\sum_{i=1}^{4}(d_i-2)\le 4$.
\proofend
\end{enumerate}
\end{lem}

\proof
(1) Note that the conclusion holds for
$(d_1,d_2,d_3)=(3,11,11)$.

Suppose that $\sum_{i=1}^{2}(d_i-2)\ge 11$.
Then, $d_2\ge 8$.
Let $d_1=3+\alpha_1$ and $d_2=8+\alpha_2$, where $\alpha_1,\alpha_2\ge 0$.
As $\sum_{i=1}^{2}(d_i-2)\ge 11$, $\alpha_1+\alpha_2\ge 4$.
As $d_3\ge d_2$,
we have
$$
\frac 12 < \sum_{i=1}^{3}\frac{1}{d_i}
\le \frac{2}{8+\alpha_2}+\frac{1}{3+\alpha_1},
$$
which implies that $4\alpha_1+\alpha_2+\alpha_1\alpha_2<4$,
a contradiction with $\alpha_1+\alpha_2\ge 4$.

(2) and (3) can be proved similarly.
\proofend

\begin{pro}\relabel{pro1}
For any planar graph $G$ with $\delta(G)\ge 3$,
there exists a cycle $C$ in $G$ with $\mu(C)\le 10$.
\end{pro}

\proof For any face $f$ of a drawing of $G$,
the weight of $f$, denoted by $\omega(f)$,  is defined as follows:
$$
\omega(f)=\sum_{v\thicksim f} \frac{1}{d_G(v)},
$$
where $v\thicksim f$ means that $v$ is a vertex on the boundary of $f$.
We have $\omega(f)\le \ell(f)/3$, where $\ell(f)$ denotes
the length of the boundary of $f$. Obviously, $|V(G)|=\sum_{f\in F(G)}\omega(f)$ and $2|E(G)|=\sum_{f\in F(G)}\ell(f)$.


By Euler's polyhedron formula,
$$
\sum_{f\in F(G)}(\omega(f)-\ell(f)/2+1)=2,
$$
implying that $\omega(f_0)-\ell(f_0)/2+1> 0$ holds for some face $f_0$.
As $\omega(f_0)\le \ell(f_0)/3$,
$$
0<\omega(f_0)-\ell(f_0)/2+1\le \ell(f_0)/3-\ell(f_0)/2+1=1-\ell(f_0)/6,
$$
implying that $\ell(f_0)\le 5$.
Let $C$ denote the boundary of face $f_0$.
Since $\delta(G)\ge 3$,  $C$ must be a cycle of $G$.

First consider the case that $\ell(f_0)=3$.
The inequality $\omega(f_0)-\ell(f_0)/2+1> 0$
implies that $\omega(f_0)>1/2$,
i.e., $1/d_1+1/d_2+1/d_3>1/2$,
where $d_1,d_2,d_3$ denote the degrees of the three vertices
on cycle $C$ with $d_1\le d_2\le d_3$.
By Lemma~\ref{lem1} (1),
$\mu(C)=d_1-2+d_2-2\le 10$.

In the case $\ell(f_0)\in \{4,5\}$,
$\omega(f_0)-\ell(f_0)/2+1> 0$ implies that $\omega(f)>\ell(f_0)/2-1$.
By Lemma~\ref{lem1}, 
there exists $v\in V(C)$ such that
$\sum_{u\in V(C)\setminus \{v\}}(d_G(u)-2)\le 5$.
Thus $\mu(C)\le 5$.
\proofend

\begin{pro}\relabel{pro2}
For any graph $G$ with
$\delta(G)\ge 3$,
$G$ contains a cycle $C$ with $\mu(C)\le sk(G)+10$.
\end{pro}

\proof Let $t=sk(G)$ and
$E_0$ be a subset of $E(G)$ with $|E_0|=t$
such that  $G\setminus E_0$ is planar.
The proof is by induction on $t$.
If $t=0$, then $G$ is planar, and the result follows from
Proposition~\ref{pro1}. 

We now assume that $t\ge 1$.
Let $e=v_1v_2\in E_0$ and $G'=G\setminus e$.
Then $sk(G')=t-1$.
For $i=1,2$, if $d_G(v_i)=3$,
then $d_{G'}(v_i)=2$ and there exists an edge $e_i$ in $G'$
incident with $v_i$ and some vertex $u_i$.
Note that $u_1$ and $u_2$ may be not distinct.

\begin{figure}[htbp]
\centering
\includegraphics[width=16 cm]{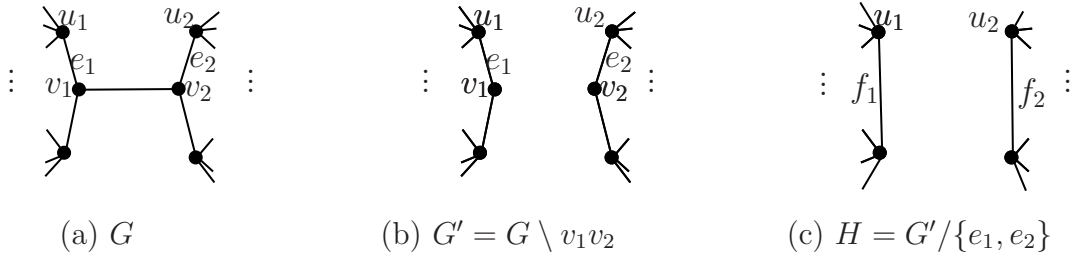}

(a) $G$  \hspace{3 cm}
(b) $G'=G\setminus v_1v_2$ \hspace{2 cm}
(c) $H=G'\slash \{e_1,e_2\}$

\caption{Graphs $G, G'$ and $H$ when $d_G(v_1)=d_G(v_2)=3$}\label{f5}
\end{figure}

Let $E_1$ be the subset of $\{e_1,e_2\}$
such that $e_i\in E_1$ \iff $d_G(v_i)=3$.
Clearly,
$|E_1|=2$ \iff $d_G(v_1)=d_G(v_2)=3$, and
$E_1=\emptyset$ \iff $d_G(v_1)>3$ and $d_G(v_2)>3$.
Let $H=G'\slash E_1$, i.e., $H$ is the graph obtained from
$G'$ by contracting all edges in $E_1$.
Thus, if $E_1=\emptyset$, then $H$ is $G'$ itself;
if $E_1=\{e_1,e_2\}$, then $H$ is as shown in Figure~\ref{f5} (c).

As $H=G'\slash E_1$ and each edge $e_i$ in $E_1$
is incident with a vertex of degree $2$ in $G'$,
we have $sk(H)=sk(G')=t-1$.
By the constriction of $H$,
$$
\delta(H)=\min\{d_{G'}(u): u\in V(G'), d_{G'}(u)>2\}\ge 3.
$$
By the inductive assumption,
$H$ contains a cycle $C_1$ with $\mu_{H}(C_1)\le t-1+10=t+9$.
Assume that $v$ is a vertex in $C_1$ such that
$$
\mu_{H}(C_1)=\sum_{u\in V(C_1)\setminus \{v\}}(d_H(u)-2).
$$
As $v$ must be a vertex in $C_1$ with the maximum degree in $H$,
we have $d_H(v)\ge 3$.

Let $C$ be the cycle of $G'$ obtained from $C_1$ by replacing edge $f_i$ by
path $P_i$ whenever $f_i$ is an edge in $C_1$ for $i=1,2$,
where $P_i$ is the path in $G'$ of length 2 which has $v_i$
as its center vertex in the case $d_{G'}(v_i)=2$.

Note that $C$ is also a cycle in $G$.
Recall that $G'=G\setminus e$ and $e$ joins $v_1$ and $v_2$.
There are two cases to consider.

\noindent {\bf Case 1}: $\{v_1,v_2\}\not\subseteq  V(C)$.

In this case,
$$
\mu_G(C)\le 1+\mu_H(C_1)\le t+10.
$$

\noindent {\bf Case 2}: $\{v_1,v_2\}\subseteq  V(C)$.

\begin{figure}[htbp]
\centering
\includegraphics[width=6 cm]{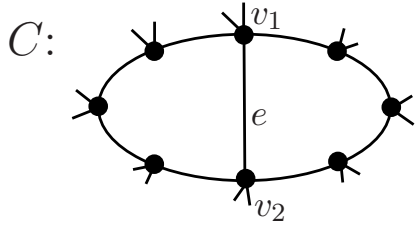}
\caption{Edge $e$ is a chord of cycle $C$}\label{f4}
\end{figure}

Note that $C\cup e$ is a subgraph in $G$, as
shown in Figure~\ref{f4}.
Let $C_1$ and $C_2$ be the two cycles in $C\cup e$
with $e\in E(C_i)$ for $i=1,2$ and $v\in V(C_1)$.

Observe that
\begin{eqnarray}\label{eq2-2}
\mu_G(C_1)
&\le& \sum_{u\in V(C_1)\setminus \{v\}}(d_{G}(u)-2)
\nonumber \\
&\le &\sum_{u\in V(C)\setminus \{v\}}(d_{G}(u)-2)
-1-(|V(C_2)|-2)(\delta(G)-2)
\nonumber \\
&=&\mu_{G'}(C)+2-1-(|V(C_2)|-2)(\delta(G)-2)
\nonumber \\
&\le &\mu_{H}(C_1)+2-1-1
\nonumber \\
&=&t+9.
\end{eqnarray}
Thus, the result holds.
 \proofend

\resection{Upper bound for $cr(G)$ in terms of $sk(G)$
\relabel{sec3}}

It is readily checked that $sk(G)\leq cr(G)$. In this section, we shall establish a new relationship between the crossing number
and skewness of a graph.

\begin{lem}\label{le02}
Let $G$ be a graph with $n\ge 4$ vertices,
and $e$ be an edge in $G$
such that $G\setminus e$ is planar.
For any planar drawing $D_1$ of $G\setminus e$,
$D_1$ can be extended to a drawing $D$ of $G$ with
$cr(D)\le \frac {2n-7}{3}$.
\end{lem}

\proof It suffices to prove this result for the case that
$G\setminus e$ is a maximal planar graph.

Let $G_1$ denote $G\setminus e$
and $D_1$ be a planar drawing of $G_1$.
As $G_1$ is a maximal planar graph, $G_1$ is $3$-connected.
Let $D_1^*$ denote the dual of $D_1$.
Thus, $D_1^*$ is 3-connected~\cite[Exercise 10.2.7]{B}.

Assume that $e$ is an edge of $G$ joining vertices $v_1$ and $v_2$ in $G$.
As $G_1$ is $3$-connected, $\delta(G_1)\ge 3$.
Thus, each vertex $v_i$ is incident with at least three faces of $D_1$.
Assume that $f_1,f_2,f_3$ are faces of $D_1$ which are
incident with $v_1$
and $f'_1,f'_2,f'_3$ are faces of $D_1$ which are
incident with $v_2$.

If $\{f_1,f_2,f_3\}\cap \{f'_1,f'_2,f'_3\}\ne \emptyset$,
say $f_1=f'_1$,
then $D_1$ can be extended to a planar drawing $D$ of $G$
by adding an edge within face $f_1$
joining $v_1$ and $v_2$ in $D_1$.
Thus $cr(D)=0\le (2n-7)/3$ holds.

Now assume that
$\{f_1,f_2,f_3\}\cap \{f'_1,f'_2,f'_3\}= \emptyset$.
As $D^*_1$ is $3$-connected,
there exist $3$ vertex-disjoint paths $P_1,P_2$ and $P_3$
in $D_1^*$ connecting vertices
in $\{f_1,f_2,f_3\}$ to vertices in $\{f'_1,f'_2,f'_3\}$.
Observe that
\begin{equation}\label{eq3-1}
\sum_{i=1}^3 |E(P_i)|=\sum_{i=1}^3 (|V(P_i)|-1)
\le |V(D_1^*)|-3=2n-4-3=2n-7,
\end{equation}
where $|V(D_1^*)|=2n-4$ follows from the fact that
$D_1$ is a maximal plane graph
and $|V(D_1^*)|$ is equal to the number of faces of $D_1$.

Assume that $|E(P_1)|\le |E(P_i)|$ for $i=2,3$.
By (\ref{eq3-1}), $|E(P_1)|\le (2n-7)/3$.
Assume that $P_1$ is a path in $D_1^*$ joining
vertices $f_1$ and $f'_1$.
As $f_1$ is a face of $D_1$ incident with $v_1$
and $f'_1$ is a face of $D_1$ incident with $v_2$,
$P_1$ actually generates a curve on the plane
connecting $v_1$ and $v_2$
which crosses with exactly $|E(P_1)|$ edges in $D_1$.
This curve represents a way of drawing edge $e$ in $D_1$.
Thus, we get a drawing $D$ of $G$
with $cr(D)=|E(P_1)|\le (2n-7)/3$.
\proofend

According to the Lemma~\ref{le02}, we now reveal the relationship between $cr(G)$ and $sk(G)$.

\begin{theo}\label{th4}
Let $G$ be a graph with $n$ vertices. Then,
\begin{equation*}\label{eq3-2}
cr(G)\le \frac{3 sk(G)^2+(4n-17)sk(G)}{6}.
\end{equation*}
Moreover, the upper bound is tight.
\end{theo}

\proof If $sk(G)=0$, then $G$ is planar and $cr(G)=0$.

Now assume that $sk(G)=t>0$.
By definition, there exists a set $E_0$ of edges in $G$
with $|E_0|=t$ such that $G\setminus E_0$ is planar.

Let $G_1$ denote the subgraph $G\setminus E_0$
and $D_1$ be a planar drawing of $G_1$.
Applying Lemma~\ref{le02} to each edge in $E_0$,
we get a  drawing $D$ of $G$ such that
$$
cr(D)\le t\frac{2n-7}3+{t\choose 2}.
$$
As $t=sk(G)$ and $cr(G)\le cr(D)$,
the claim follows.
It can be verified easily that the upper bound is tight for the complete graph $K_{5}$.
\proofend

\resection{Main results on crossing-critial graphs\relabel{sec4}
}

Recall that $\M_k$ denotes the set of $k$-crossing-critical graphs.

\begin{theo}\relabel{th5}
Let $G$ be a graph in $\M_k$ with minimum degree $\delta$.
If $G$ contains a cycle $C$ with $\mu(C)=s$,  then,
$$ cr(G)\le\left\{
\begin{array}{ll}
\displaystyle
2k+\frac{s-5}{2}, &\mbox{if }\delta=3; \\
\displaystyle
2k-sk(G)+\frac{\delta(s-\delta+2)}{2(\delta-2)}, \quad &\mbox{if }
\delta\ge 4.
\end{array}
\right.
$$
\end{theo}

\proof
As $\mu(C)=s$, 
there exists  some vertex $v$ in $C$ such that
$$
|E(P)|=|V(C)|-1\le
\frac 1{\delta-2}\sum_{u\in V(C)\setminus \{v\}}(d_G(u)-2)
= \frac s{\delta-2}.
$$
Let $e$ be an edge of $C$ with ends $v$ and $w$
and $P$ denote the path of $C\setminus e$. As $G$ is $k$-crossing-critical, $cr(G\setminus e)\le k-1$.
Let $D$ be a drawing of $G\setminus e$ with at most $k-1$ crossings.
Note that edges in $P$ may cross each other in the drawing $D$.
We regard the drawing of $P$ as a planar graph $H$ with vertices of degrees
2 and 4. Let $P'$ be a shortest path in $H$ joining $v$ and $w$. There are two ways of reconnecting $v$ and $w$ close to $P'$, one for each side of $P'$.

Let $r_D(P)$ denotes the number of crossings of edges of $P$ in $D$. It is not hard to verify that the total number of crossings in these two drawings
of $e$ is at most
\begin{eqnarray*}
  \lambda &=& \sum_{u\in V(C)\setminus\{v,w\}}(d_G(u)-2)+2r_D(P)\\
   &\le& s-(\delta-2)+2r_D(P)
\end{eqnarray*}
Therefore, one of the two drawings for $e$ crosses
at most $\lambda/2$ edges of $D$.
Note that $r_D(P)\le k-1$,
implying that
\begin{equation}\label{eq2}
cr(G)\le \frac{\lambda}{2}+k-1
\le \frac{s-\delta+2}{2}+2k-2
=\frac{s-\delta}{2}+2k-1.
\end{equation}

If $\delta=3$, then it follows that
$cr(G)\le 2k+(s-5)/2$. 

Now we consider the case that $\delta\ge 4$. Removing from $D$ the edges (at most $s/(\delta-2)$) of $P$
leaves a drawing with at most $k-1-r_D(P)$ crossings. Therefore, there is a set
of at most $1+s/(\delta-2)+k-1-r_D(P)$ edges whose removal from $G$ leaves a
planar graph. Thus,$sk(G)\le 1+s/(\delta-2)+k-1-r_D(P)$,
implying that
$$r_D(P)\le \frac s{\delta-2}+k-sk(G).$$
As $\lambda \le s-(\delta-2)+2r_D(P)$, by (\ref{eq2}),
\begin{eqnarray}
cr(G)&\le &\frac{\lambda}{2}+k-1
\nonumber \\
&\le& \frac{s-\delta+2}{2}+r_D(P)+k-1
\nonumber \\
&\le &\frac{s-\delta+2}{2}+k-1 +\frac s{\delta-2}+k-sk(G).
\end{eqnarray}
The result holds.
\proofend

\begin{theo}\label{th3}
Let $G\in \M_k$ with $n$ vertices and minimum degree $\delta$. Then,
$$
cr(G)\le\left\{
\begin{array}{lll}
\displaystyle
2.5(k+1), &&\mbox{if }\delta=3;\\
2(k+4), &  &\mbox{if }\delta= 4;\\
\displaystyle
2k-\sqrt k/2n+35/6, &  &\mbox{if } \delta\ge 5.
\end{array}
\right.
$$
\end{theo}

\proof Let $t=sk(G)$. By Proposition~\ref{pro2},
$G$ contains a cycle $C$ with $\mu(C)\le t+10$.

If $\delta=3$, by Theorem~\ref{th5},
$cr(G)\le 2k+(t+10-5)/2=2k+2.5+0.5t\le 2.5k+2.5$,
as $t\le k$.

If $\delta=4$, by Theorem~\ref{th5},
$cr(G)\le 2k-t+(t+10-4+2)=2k+8.$

Now consider the case that $\delta\ge 5$.
by Theorem~\ref{th5},
\begin{equation}\label{eq3}
cr(G)\le 2k-t+\frac{\delta(t+10-\delta+2)}{2(\delta-2)}
\le 2k-t+\frac{5(t+7)}{6}
=2k+\frac{35-t}{6}.
\end{equation}
By (\ref{eq3}), if $t\ge \frac {3\sqrt k}n$, the result holds.
In the following, assume that  $t<\frac {3\sqrt k}n$.

By Theorem~\ref{th4},
\begin{equation}\label{eq4}
cr(G)\le \frac{3t^2+(4n-17)t}{6}
<\frac{3(9k/n^2)+(4n-17)\frac {3\sqrt k}n}{6}
=\frac{9k+(4n-17)n\sqrt k}{2n^2}.
\end{equation}
If $k=1$, then, by (\ref{eq4}),
\begin{equation}\label{eq5}
cr(G)\le \frac{9+(4n-17)n}{2n^2}<2,
\end{equation}
and the result holds.

If $k\ge 2$, then, by (\ref{eq4}),
\begin{eqnarray}\label{eq6}
& &cr(G)-\left (2k+35/6-\frac {1}{2n} \sqrt k \right)
\nonumber \\
&\le &\frac{9k+(4n-17)n\sqrt k}{2n^2}
-\left (2k+35/6-\frac {1}{2n} \sqrt k \right)
\nonumber \\
&<&\frac{9k+(4n-0)n\sqrt k}{2n^2}
-\left (2k+0-\frac {1}{2n} \sqrt k \right)
\nonumber \\
&=&\frac{(4.5-2n^2)k+(2n+0.5)n\sqrt k}{n^2}
\nonumber \\
&<&0,
\end{eqnarray}
where the last inequality follows from the facts
that the solution of the inequality
$(4.5-2n^2)k+(2n+0.5)n\sqrt k<0$ is $k>n^2(4n+1)^2/(4n^2-9)^2$
and that $2>n^2(4n+1)^2/(4n^2-9)^2$ holds for all $n\ge 5$.
Thus, we complete the proof.
\proofend

\end{document}